\documentclass{agtart_a}
\pdfoutput=1

\title[The diameter of the set of boundary slopes of a knot]{The diameter of the set of boundary slopes\\of a knot}

\author{Ben Klaff}
\givenname{Ben}
\surname{Klaff}
\address{Department of Mathematics\\
University of Texas at Austin\\\newline
1 University Station\\
Austin, TX 78741\\USA}
\email{klaff@math.utexas.edu}
\urladdr{}

\author{Peter B Shalen}
\givenname{Peter B}
\surname{Shalen}
\address{Department of Mathematics, Statistics, and 
Computer Science (M/C 249)\\\newline
University of Illinois at Chicago\\
851 S. Morgan St.\\
Chicago, IL 60607-7045\\USA}
\email{shalen@math.uic.edu}
\urladdr{}

\volumenumber{6}
\issuenumber{}
\publicationyear{2006}
\papernumber{40}
\startpage{1095}
\endpage{1112}

\doi{}
\MR{}
\Zbl{}


\received{12 November 2005}
\revised{}
\accepted{14 March 2006}
\published{}
\publishedonline{29 August 2006}
\proposed{}
\seconded{}
\corresponding{}
\editor{Martin Scharlemann}
\version{}

\arxivreference{math.GT/0412147}

\makeatletter
\def\cnewtheorem#1[#2]#3{\newtheorem{#1}{#3}[section]
\expandafter\let\csname c@#1\endcsname\c@para}
\def\dnewtheorem#1[#2]#3{\newtheorem{#1}{#3}
\expandafter\let\csname c@#1\endcsname\c@thm}
\makeatother

\AtBeginDocument{\let\bar\wbar\let\tilde\wtilde\let\hat\what}
\makeop{frontier}

\theoremstyle{plain}
\newtheorem{thm}{Theorem}

\dnewtheorem{cor}[thm]{Corollary}

\swapnumbers
\theoremstyle{remark}
\newtheorem{para}{\hspace{-4pt}}[section]
\cnewtheorem{remark}[para]{Remark}
\cnewtheorem{remarks}[para]{Remarks}
\cnewtheorem{definition}[para]{Definition}
\cnewtheorem{definitions}[para]{Definitions}
\cnewtheorem{notation}[para]{Notation}

\makeautorefname{para}{}

\theoremstyle{plain}
\cnewtheorem{theorem}[para]{Theorem}
\cnewtheorem{lemma}[para]{Lemma}
\cnewtheorem{proposition}[para]{Proposition}
\cnewtheorem{corollary}[para]{Corollary}

\numberwithin{equation}{para}

\newcommand\Number{\begin{para}}
\newcommand\EndNumber{\end{para}}
\newcommand\Remark{\begin{remark}}
\newcommand\EndRemark{\end{remark}}
\newcommand\Lemma{\begin{lemma}}
\newcommand\EndLemma{\end{lemma}}
\newcommand\Proposition{\begin{proposition}}
\newcommand\EndProposition{\end{proposition}}
\newcommand\Proof{\begin{proof}}
\newcommand\EndProof{\end{proof}}

\newcommand\inter{\mathop{\rm int}}
\newcommand\qinfty{\Q\cup\{\infty\}}
\newcommand\sea{c}

\begin{document}

\begin{asciiabstract}
Let K be a tame knot with irreducible exterior M(K) in a closed,
connected, orientable 3--manifold Sigma such that pi_1(Sigma)
is cyclic. If infinity is not a strict  boundary slope, then the
diameter of the set of strict boundary slopes of K, denoted d_K, is
a numerical invariant of K.  We show that either (i) d_K >= 2 or (ii)
K is a generalized iterated torus knot. The proof combines results from
Culler and Shalen [Comment. Math. Helv. 74 (1999) 530-547] with a result
about the effect of cabling on boundary slopes.
\end{asciiabstract}

\begin{htmlabstract}
Let K be a tame knot with irreducible exterior M(K) in a
closed, connected, orientable 3&ndash;manifold &Sigma; such that
&pi;<sub>1</sub>(&Sigma;) is cyclic. If &infin; is not a strict
boundary slope, then the diameter of the set of strict boundary
slopes of K, denoted d<sub>K</sub>, is a numerical invariant of K.
We show that either (i) d<sub>K</sub>&ge;2 or (ii) K is a generalized
iterated torus knot. The proof combines results from Culler and Shalen
[Comment. Math. Helv. 74 (1999) 530&ndash;547] with a result about the
effect of cabling on boundary slopes.
\end{htmlabstract}

\begin{webabstract}
Let $K$ be a tame knot with irreducible exterior $M(`K)$ in a closed,
connected, orientable 3--manifold $\Sigma$ such that $\pi_{1}(\Sigma)$ is
cyclic. If $\infty$ is not a strict  boundary slope, then the diameter of
the set of strict boundary slopes of $K$, denoted $d_{K}$, is a numerical
invariant of $K$.  We show that either (i) $d_{K} \geq 2$ or (ii) $K$
is a generalized iterated torus knot. The proof combines results from
Culler and Shalen [Comment. Math. Helv. 74 (1999) 530-547] with a result
about the effect of cabling on boundary slopes.
\end{webabstract}

\begin{abstract}
Let $K$ be a tame knot with irreducible exterior $M(`K)$ in a closed,
connected, orientable 3--manifold $\Sigma$ such that $\pi_{1}(\Sigma)$ is
cyclic. If $\infty$ is not a strict  boundary slope, then the diameter of
the set of strict boundary slopes of $K$, denoted $d_{K}$, is a numerical
invariant of $K$.  We show that either (i) $d_{K} \geq 2$ or (ii) $K$
is a generalized iterated torus knot. The proof combines results from
Culler and Shalen \cite{slopediff} with a result about the effect of
cabling on boundary slopes.
\end{abstract}

\maketitle

\section*{Introduction}

Let $K$ be a (tame) knot in a connected, closed, orientable $3$--manifold
$\Sigma$, such that the exterior $M(`K)$ of $K$ is irreducible and
$M(`K)$ contains no strict essential surface with meridian boundary
slope.  The diameter $d_{K}$ of the set of all boundary slopes of
strict essential surfaces in $M(`K)$ is a natural invariant of $K$.
(The definition of $d(`K)$ and of other specialized terms used here
will be reviewed below.)  Properties of $d_{K}$ have topological and
algebraic meaning for the knot $K$.  For example, the main result of
Culler and Shalen \cite{oldcs}, which implies Neuwirth's conjecture \cite{Neuwirth} that
classical knot groups are nontrivial amalgamated free products, is
equivalent to the assertion that under suitable mild restrictions on
$K$, we have $d_{K} \neq 0$.  (For further discussion of this
connection, see Shalen \cite{handbook}.)

It was shown in Culler and Shalen \cite{slopediff} that if $\pi_1(\Sigma)$ is trivial
and the knot $K$ is nontrivial, then $d_K\ge2$.  In this paper we
extend this result to the case in which $\pi_1(\Sigma)$ is cyclic.
\fullref{theoremB} below asserts that in this situation we still have
$d_K\ge2$, unless the knot $K$ belongs to a certain special class of
knots that we call generalized iterated torus knots.

When $K$ belongs to the special class just mentioned, then a great
deal is known about $K$ and also $\Sigma$: it is not hard to classify
all generalized iterated torus knots, to show that a manifold
containing such a knot must be a lens space, and to calculate $d_{K}$
when $K$ is such a knot.  In many (but not all) cases, it turns out
that $d_{K} < 2$.

\fullref{theoremB} is proved by combining one of the results of
\cite{slopediff} with another new result, \fullref{theoremA}, which
asserts that a $q$--strand cabling of a knot increases the invariant
$d_K$ by a factor of at least $q^2$. 

By combining \fullref{theoremA} and \fullref{theoremB} we obtain \fullref{corC} which gives
a stronger version of \fullref{theoremB} for the case of a cable knot.
For example, it follows from \fullref{corC} that if $K$ is a cable of a
hyperbolic knot in $S^3$ then $d_K\ge8$.

\fullref{corC} will be used in Klaff \cite{thesis}, where it is shown that
if $\pi_1(\Sigma)$ is a finite cyclic group of odd order, and if the
knot $K\subset\Sigma$ is not a generalized iterated torus knot, then
$d_K$ is strictly greater than $2$.

One of the ingredients in the proof of \fullref{theoremB} is \fullref{mfldwithtorus}, which provides a criterion for deciding whether a
surface in a cabled knot exterior---or more generally, in an
irreducible, orientable $3$--manifold that contains an essential
torus---is strict and essential. This result is of independent
interest, and is applied by Culler and Shalen \cite{twosurfs} in a rather different
context.

Before giving formal statements
of Theorems A and B, we shall set up some general conventions and
define some of the terms we used above.

We shall work in the piecewise linear category throughout this paper.

\Number\label{unislope}
If $T$ is a $2$--dimensional torus, we define a \textit{slope\/} on $T$ to
be an isotopy class of homotopically nontrivial simple closed curves
in $T$. The set of all slopes on $T$ will be denoted by ${\mathcal S}(T)$.
  
The isotopy classes of homotopically nontrivial \textit{oriented\/} simple
closed curves in $T$ are in natural bijective correspondence with
elements of $H_1(T;\Z)$ which are \textit{primitive\/} in the sense of
not being divisible by any integer greater than $1$. Thus there is a natural
two-to-one map from the set of primitive elements of $H_1(T;\Z)$
onto  ${\mathcal S}(T)$. We shall denote this map by
$\alpha\mapsto\langle\alpha\rangle$. We have
$\langle\alpha\rangle=\langle\alpha'\rangle$ if and only if
$\alpha'=\pm\alpha$.
\EndNumber
 
\Number\label{multislope} If $C$ is a nonempty closed $1$--manifold in
a $2$--torus $T$, and $C$ has no homotopically trivial components, then
all components of $C$ have the same slope $\sigma\in{\mathcal S}(T)$. We
call $\sigma$ the \textit{slope\/} of $C$.  
\EndNumber

\Number\label{essential}
A $3$--manifold $M$ is \textit{irreducible\/} if $M$ is connected and every
$2$--sphere in $M$ bounds a ball. 

An \textit{essential surface\/} in an irreducible, orientable $3$--manifold
$M$ is a two-sided properly embedded surface in $M$ which is nonempty
and $\pi_1$--injective, and has no $2$--sphere components and no
boundary-parallel components.
\EndNumber

\Number\label{strict}
Suppose that $M$ is a compact, orientable irreducible $3$--manifold
whose boundary components are tori.  A connected essential surface in
$M$ is called a \textit{semifiber\/} if either $F$ is a fiber in a
fibration of $M$ over $S^1$, or $F$ is the common frontier of two
$3$--dimensional submanifolds of $M$, each of which is a twisted
$[0,1]$--bundle with associated $\{0,1\}$--bundle $F$. An essential
surface $F\subset M$ is termed \textit{strict\/} if no component of $F$ is
a semifiber.  A strict essential surface has no disk components,
since an irreducible knot manifold which has an essential disk must be
a solid torus, and the essential disk in a solid torus is a fiber.
\EndNumber

\Number\label{meetsall}
Since a semifiber in a bounded $3$--manifold $M$ must meet every
component of $\partial M$, any essential surface that is disjoint from
at least one component of $\partial M$ must be strict.
\EndNumber

\Number \label{boundaryslope} Let $M$ be a compact orientable
$3$--manifold, and let $T$ be a component of $\partial M$ which is a
torus.  If $F$ is an essential surface in $M$ that meets $T$, then
$\partial F \cap T$ is a $1$--manifold in $T$ having no homotopically
trivial components. Thus by \ref{multislope}, $\partial F \cap T$ has
a well-defined slope $\sigma\in{\mathcal S}(T)$, which we call the \textit{boundary slope of $F$ on $T$\/}.  
\EndNumber

\Number\label{Hatcher} We define a \textit{knot manifold\/} to be a
connected, compact, orientable $3$--manifold $M$ such that $\partial M$
is a torus.  If $M$ is a knot manifold, we define a \emph{(strict)
  boundary slope of $M$} to be an element of ${\mathcal S}(\partial M)$
which arises as the boundary slope of some bounded (strict) essential
surface $F$ in $M$.  A theorem of Hatcher's \cite{Hatcher,JSe} implies that
for any given knot manifold $M$, the boundary slopes of $M$ form a
finite subset of ${\mathcal S}(\partial M)$. In particular, the strict boundary
slopes of $M$ form a finite subset of ${\mathcal S}(\partial M)$.  \EndNumber

\Number\label{meridian} If $K$ is a (PL) knot in a closed, orientable
$3$--manifold $\Sigma$, we shall denote by $V(`K)$ a regular
neighborhood of $K$, and by $M(`K)$ the \textit{exterior\/} of $K$, defined
by $M=\overline{\Sigma-V(`K)}$.  Note that $M(`K)$ is a knot
manifold. Since $V(`K)$ and hence $M(`K)$ are
well-defined up to ambient isotopy in $\Sigma$, and in particular up
to homeomorphism, our main results are independent of the choice of
$V(`K)$. In general we shall implicitly suppose an arbitrary choice of
$V(`K)$ to have been made, but at one point in \fullref{Asection}
it will be necessary to be more explicit.

A \textit{meridian\/} of $K$ is a nontrivial simple closed curve in the
torus $\partial M(`K)$ which bounds a disk in $V(`K)$.  Such a curve
exists and is unique up to isotopy.  Thus there is a well-defined \textit{meridian slope\/} in $\partial M$. A primitive element $\mu$ of
$H_1(\partial M(`K);\Z)$ is called a \textit{meridian class\/} for $K$
if $\langle\mu\rangle$ is the meridian slope. Thus $K$ has exactly two
meridian classes, and they differ by a sign.

\Number
The knot $K$ will be termed \textit{meridionally small\/} if its meridian
slope is not a boundary slope for $M(`K)$. (In the case where $\Sigma$
is an irreducible nonHaken manifold, meridionally small knots in
$\Sigma$ are sometimes called ``smallish knots.'')
\EndNumber

\Number\label{numerical}
We define a \textit{framing\/} for $K$ to be an ordered basis
$(\mu,\lambda)$ for $H_1(\partial M(`K);\Z)$ such that $\mu$ is a
meridian class.  
\EndNumber

If $(\mu,\lambda)$ an arbitrary framing for $K$, there is a
bijection $\nu=\nu_{\mu,\lambda}$ from ${\mathcal S}(\partial M(`K))$ to $\Q 
\cup\{\infty\}$ defined by setting
$\nu(\langle\alpha\rangle)=\omega(\alpha,\lambda)/\omega(\alpha,\mu)$
where $\alpha\in H_1(\partial M(`K); \Z)$ is any primitive
and $\omega$ denotes homological intersection number. Equivalently,
if we write $\alpha=a\mu+b\lambda$, where $a$ and $b$ are relatively
prime integers, then $\nu(\langle\alpha\rangle)=-a/b$.
\EndNumber

\Number\label{translation} Note that if $(\mu_1,\lambda_1)$ and
$(\mu_2,\lambda_2)$ are two framings for the knot $K$, then there
exist $h\in\Z$ and $\epsilon\in\{\pm1\}$ such that
$\mu_2=\epsilon\mu_1$ and $\lambda_2=\lambda_1+h\mu_1$. It follows
that if $\sigma$ is any slope on $\partial M(`K)$, and if we set
$s_i=\nu_{\mu_i,\lambda_i}(\sigma)$ for $i=1,2$, then we have
$$s_2=\epsilon s_1+h.$$
\EndNumber

\Number\label{NumHatch} Now suppose that $M(`K)$ is irreducible. If $F$
is a bounded essential surface in $ M(`K)$, we define the \textit{  numerical boundary slope\/} of $F$, with respect to a given framing
$(\mu,\lambda)$, to be the image of the boundary slope of $F$ under
the bijection $\nu_{\mu,\lambda}\co {\mathcal S}(\partial M(`K)) \to \Q 
\cup\{\infty\}$.  In analogy with \ref{Hatcher}, we define a
\emph{(strict) numerical boundary slope} of $K$ to be a slope on
$\partial M(`K)$ which arises as the boundary slope of some bounded
(strict) essential surface $F$ in $M$.  If we denote by ${\mathcal
  B}_{\mu,\lambda}(`K)\subset\qinfty$ the set of all numerical boundary
slopes of bounded strict essential surfaces in $M(`K)$ with respect to
the framing $(\mu,\lambda)$, the theorem of Hatcher's quoted in
\ref{Hatcher} implies that ${\mathcal B}_{\mu,\lambda}(`K)$ is a finite
set.

In particular, if $K$ is meridionally small and has at least one strict
boundary slope, then ${\mathcal B}_{\mu,\lambda}(`K)\subset\Q$ is a finite,
nonempty subset of $Q$; thus in this case ${\mathcal
  B}_{\mu,\lambda}(`K)\subset\Q$ has a
greatest element $s_{\rm max}(`K,\mu,\lambda)$ and a least element
$s_{\rm min}(`K,\mu,\lambda)$.

The observation \ref{translation} shows that if $(\mu_1,\lambda_1)$
and $(\mu_2,\lambda_2)$ are framings of $K$, the sets ${\mathcal
  B}_{\mu_1,\lambda_1}(`K)$ and ${\mathcal B}_{\mu_2,\lambda_2}(`K)$
differ only by an integer translation and a possible change of sign.
In particular, in the case where $K$ is meridionally small and there
is at least one strict boundary slope, the \textit{diameter\/} $d=s_{\rm
  max}(`K,\mu,\lambda)-s_{\rm min}(`K,\mu,\lambda)$ of $ {\mathcal
  B}_{\mu,\lambda}(`K)$ is independent of the framing $(\mu,\lambda)$
and is therefore an invariant of the knot $K$, which we denote by
$d_K\in\Q$.  If there are no strict boundary slopes for $K$ we
set $d_K=-\infty.$ Thus the invariant $d_K$ is defined for every
meridionally small knot $K$ in a closed, orientable, irreducible
$3$--manifold $\Sigma$.  \EndNumber

\Number\label{round}
A knot $K\subset\Sigma$ is said to be
\emph{round} if some $M(`K)$ admits a solid torus as a connected
summand. Thus when $M(`K)$ is assumed to be irreducible, $K$ is round
if and only if $M(`K)$ is a solid torus. Note that this implies that
$\Sigma$ has a genus--$1$ Heegaard splitting; that is, $\Sigma$ is
homeomorphic either to $S^2\times S^1$ or to a (possibly trivial) lens
space.
\EndNumber

\Number\label{cable} 
A knot $K'\subset \Sigma$ is called a \emph{cabling} of a knot
$K\subset\Sigma$ if there exists a regular neighborhood $U$ of $K$ such
that $K'$ is a simple closed curve on the boundary of $U$, and the
geometric intersection number $q$ of $K'$ with the boundary of a
meridian disk for the solid torus $U$ is greater
than or equal to $2$. We shall refer to the integer $q\ge2$ as the
\textit{number of strands\/} of the cable $K'$.

Note that according our definition, a knot $L$ which is isotopic in
$\Sigma$ to a ($q$--strand) cabling of $K$ need not itself be a cabling
of $K$. However, such a knot $L$ is clearly a ($q$--strand) cabling of
some knot isotopic to $K$.

We define a \textit{$q$--strand cable knot\/} in $\Sigma$ to be a knot
$K\subset\Sigma$ such that (a) $K$ is not round and (b) $K$ is a
$q$--strand cabling of some knot in $\Sigma$. (Note that (b) does not
imply (a), since a trivial knot in $S^3$ is a $q$--strand cabling of
another trivial knot for any $q\ge2$.) We call $K$ a \textit{cable knot\/}
if it is a $q$--strand cable knot for some $q\ge2$.  
\EndNumber

\Number\label{gitk}
A knot $K\subset\Sigma$ is called a \emph{generalized iterated torus
  knot} if for some integer $n\ge0$ there exist knots $K_{0}, K_{1},
\ldots, K_{n}$ in $\Sigma$ such that (i) $K=K_{0}$, (ii) the knot
$K_{n}$ is round and (iii) for each $i$ with $0 \leq i \leq n-1$, the
knot $K_{i}$ is a cabling of $K_{i{+}1}$.  It follows from the
observation made in \ref{round} that if $\Sigma$ contains a
generalized iterated torus knot whose exterior is irreducible, then
$\Sigma$ is either a homeomorph of $S^2\times S^1$ or a (possibly
trivial) lens space.  \EndNumber

We are now in a position to give formal statements of our main results.

\begin{thm}\label{theoremA} Let $\Sigma$ be a closed,
  connected, orientable $3$--manifold, $K\subset\Sigma$ be a
  nonround knot, $q\ge2$ be an integer and
  $K'$ be a $q$--strand cabling of $K$.  Suppose that $K'$ is 
  meridionally small. Then $K$ is
  meridionally small, and $d_{K'} \geq q^2d_K$.
\end{thm}

This will be proved in \fullref{Asection}, using foundational material
that will be presented in \fullref{essentialsurfacesection}.

\begin{thm}\label{theoremB} Let $\Sigma$ be a closed,
  connected, orientable $3$--manifold such that $\pi_1(\Sigma)$ is
  cyclic. Suppose that $K\subset\Sigma$ is a meridionally small knot.
  Then either (i) $d_{K} \geq 2$ or (ii) $K$ is a generalized
  iterated torus knot.
\end{thm}

This will be proved in \fullref{Bsection} by combining
\fullref{theoremA} with results from \cite{slopediff}.

Theorems A and B taken together immediately yield the following:

\begin{cor}\label{corC} Let $\Sigma$ be a closed,
  connected, orientable $3$--manifold such that $\pi_1(\Sigma)$ is
  cyclic. Suppose that $q\ge 2$ is an integer, and that
  $K\subset\Sigma$ is a  $q$--strand cable 
  knot which is meridionally small. Then either (i)
  $d_{K} \geq 2q^2$ or (ii) $K$ is a generalized iterated torus knot.
\end{cor}

The first author was partially supported by the NSF VIGRE program and by the Chaire de recherche du Canada en alg\`ebre combinatoire et informatique math\'ematique at the Universit\'e du Qu\'ebec \`a Montr\'eal. The second author was partially supported by NSF grant DMS 0204142.

\section{Strict essential surfaces in toroidal manifolds}
\label{essentialsurfacesection}

The main result of this section, \fullref{mfldwithtorus}, will
provide a criterion for deciding whether a surface in the cabled knot
exterior is strict and essential. This is needed for the proof of
\fullref{theoremB}. The result will be proved in a more general setting: rather
than considering only cabled knot exteriors, we shall consider
arbitrary compact, irreducible, orientable $3$--manifolds that contain
essential tori. The result seems likely to be of broader interest in
$3$--manifold theory.

The following result, which will be used in the proof of \fullref{mfldwithtorus}, is also of more general interest. It says that an
essential surface in an irreducible knot manifold is
boundary-incompressible in a strong, homotopy-theoretic sense.

\Proposition\label{boundaryincompressible}
Suppose that $F$ is a bounded essential surface in an irreducible knot
manifold $M$, and suppose that $\alpha$ is a path in $F$ which has its
endpoints in $\partial F$ and is fixed-endpoint homotopic in $M$ to a
path in $\partial M$. Then $\alpha$ is fixed-endpoint homotopic in $F$
to a path in $\partial F$.
\EndProposition

\Proof
Let $C^+$ and $C^-$ denote the upper and lower semicircles
in $S^1=\partial D^2$. The hypothesis implies that there is a map
$f\co D^2 \to M$ such that $f\vert_{C^+}$ is a reparametrization of $\alpha$ and
$f(C^-)\subset\partial M$. We may choose $f$ so that $f^{-1}(\partial
M)=C^-$ and $C^+$ is a component of $f^{-1}(F)$, and so that
$f\vert_{(D^2-C^+)}$ is transverse to $F$. Among all maps with these
properties, we suppose $F$ to be chosen so as to minimize the number
of components of $f^{-1}(F)$. Since $F$ is $\pi_1$--injective, the
minimality implies that each component of $f^{-1}(F)$ is an arc. Hence
some component $A^+$ of $f^{-1}(F)$ is ``outermost'' in the sense that
$A^+$ is the frontier of a disk $\Delta\subset D^2$ with
$\Delta\cap\partial D^2\subset\inter C^-$ and $\Delta\cap
f^{-1}(F)=A^+$. A priori, $A^+$ may or may not coincide with $C^+$.

If $M'$ denotes the manifold obtained by splitting $M$ along $F$, and
$q\co M' \to M$ denotes the quotient map, we may write $f\vert_\Delta=q\circ g$
for some map $g\co \Delta \to M'$. We set $A^-=(\partial\Delta)-\inter A^+$.
We may regard $g\vert_{A^-}$ as a reparametrization of a path $\beta$ in
$q^{-1}(\partial M)$. Since $\partial M$ is a torus, the component $B$
of $q^{-1}(\partial M)$ containing $\beta(I)$ is an annulus. We claim
that $\beta$ has both its endpoints in the same component of $\partial
B$.

Suppose to the contrary that the endpoints of $\beta$ are in different
components of $\partial B$. Then $g$ is homotopic rel $A^+$ to a map
$g'\co \Delta \to M'$ such that $g'\vert_{A^-}$ is injective and $l^-=g'(A^-)$ is
a properly embedded arc in $B$ which meets some core curve of $B$
transversally in a single point. In particular $g'\vert_{\partial\Delta}$ is
a homotopically nontrivial map of $\partial\Delta$ into $\partial
M'$, and $g'(\partial\Delta)\subset l^-\cup F'$ for some component
$F'$ of $q^{-1}(F)$. It now follows from Henderson's version of
the loop theorem \cite[Theorem III.5]{henderson} 
that there is a disk $E\subset M'$
such that $\partial E$ is a nontrivial simple closed curve in
$\partial M'$ and $\partial E\subset l^-\cup F'$. We cannot have
$\partial E\subset F'$, since $F$ is $\pi_1$--injective in $M$. Hence
$\partial E$ must have the form $l^-\cup l^+$ for some properly
embedded arc $l^+$ in $F'$.

If $N$ denotes a regular neighborhood of $E$ in $M'$ then $R=N\cap F'$
is a regular neighborhood of $l^+$ in $F'$, and the closure of
$(\frontier_{\mskip1mu M'}N)-R$ is a disjoint union of disks $G_1$ and
$G_2$. The surface $F_1=(F-q(R))\,\cup\,q(G_1\cup G_2)$ is properly embedded in
$M$, and is $\pi_1$--injective since $F$ is. But since $l$ joins different
components of $\partial B$, the component of $\partial F_1$ contained
in $q(B)$ is a homotopically trivial simple closed curve in $\partial
M$. Hence $F_1$ is a disk, and by irreducibility it is the frontier of
a ball $K\subset M$. We must have either $q(N)\subset K$, in which
case $F$ is an annulus contained in the ball $K$, or $q(N)\cap \inter
K=\emptyset$, in which case $F$ is a boundary-parallel annulus. In
either case we have a contradiction to the essentiality of $F$. Thus
our claim is proved.

We may regard $f\vert_{A^+}$ as a reparametrization of a path $\alpha_0$ in
$F$. In the case where $A^+=C^+$, we may take $\alpha_0=\alpha$. Since
the endpoints of $\beta$ lie in the same component of $\partial B$,
and since $F$ is $\pi_1$--injective in $M$, the path $\alpha_0$ is
fixed-endpoint homotopic in $q(B)$ to a path $\beta_1$ in $\partial
B\subset\partial F$.  This implies the conclusion of the proposition
in the case where $A^+=C^+$. If $A^+\ne C^+$, we may use the homotopy
between $\alpha_0$ and $\beta_1$ to replace the map $f$ by a map
$f_1\co \Delta \to M$ which agrees with $f$ on $\Delta$ and maps
$D^2-\Delta$ into $F$; by perturbing $f_1$ slightly we obtain a map
$f_2$ such that $f_2^{-1}(F)$ has fewer components than $f^{-1}(F)$,
in contradiction to the minimality property of $f$. Thus the case
$A^+\ne C^+$ does not occur, and the proof is complete.
\EndProof

\Number\label{strictlyvertical} 
The proof of the main result of this section,
\fullref{mfldwithtorus}, also involves some basic facts about
incompressible surfaces in interval bundles over surfaces, which we
shall now summarize.

Suppose that $J$ is an orientable $3$--manifold which is a
$[0,1]$--bundle over a surface.  We will define a surface in $J$ to be
\textit{vertical\/} if it is a union of fibers, and to be \textit{horizontal\/}
if it is everywhere transverse to the fibers.  The \textit{vertical
  boundary} of $J$ is the inverse image of the boundary of the base of
the $I$--bundle $J$ under the projection map.

Any essential vertical annulus in $J$ is the inverse image under the
 fibration map of an essential simple closed curve in the base.

Suppose that $J$ is a trivial $[0,1]$--bundle and that $F$ is a properly
embedded $\pi_1$--injective surface in $J$ such that all components of
$\partial F$ are contained in the same component $C$ of the $\{
0,1\}$--bundle associated to $J$.  It follows from
Waldhausen \cite[Proposition 3.1]{waldhausen} that $F$ is isotopic to a
horizontal surface by an ambient isotopy that preserves the vertical
boundary of $J$, and that each component of $F$ is parallel to a
subsurface of $C$. 

As a consequence of this fact we observe that if $J$ is a trivial
$[0,1]$--bundle, and $F$ is a properly embedded $\pi_1$--injective
surface in $J$ such that $\partial F$ is contained in the vertical
boundary of $J$, then $F$ is isotopic to a horizontal surface by an
ambient isotopy that preserves the vertical boundary of $J$.

Suppose that $J$ is a $[0,1]$--bundle and that $A$ is a disjoint union of
properly embedded annuli in $J$ none of which is parallel to an
annulus contained in the $\{ 0,1\}$--bundle associated to $J$.  It
follows from \cite[Lemma 3.4]{waldhausen} in the case that $J$ is a
trivial $[0,1]$--bundle, and from \cite[Lemma 2]{brittenham} in the
twisted case that $A$ is isotopic to a vertical surface.

Suppose that $F$ is a properly embedded $\pi_1$--injective surface in
a $[0,1]
$--bundle $J$ such that $\partial F$ is contained in the vertical
boundary of $M$.  Then $F$ is isotopic to a horizontal surface.  This
follows from \cite[Proposition 3.1 and Proposition 4.1]{waldhausen}.
\EndNumber

\Proposition\label{mfldwithtorus}
Let $M$ be a compact orientable irreducible $3$--manifold containing an
essential torus $T$, let $M'$ be the manifold obtained by splitting
$M$ along $T$ and let $q\co M' \to M$ denote the quotient map.  Let $F$ be
a connected properly embedded surface in $M$ which is not isotopic to
$T$.  Then $F$ is a strict essential surface if and only if it is
isotopic to a surface $S$ transverse to $T$ such that
\begin{enumerate}
  \item each component of $q^{-1}(S)$ is essential in the component
  of $M'$ containing it;
  \item some component of $q^{-1}(S)$ is a strict essential surface
  in the component of $M'$ containing it.
  \end{enumerate}
\EndProposition

\Proof
Given a surface $S$ transverse to $T$ we will set $S'= q^{-1}(S)$ and
$T' = q^{-1}(T)$.  We let $M_1$ denote the manifold obtained by
splitting $M$ along $S$ and denote the quotient map by $p\co M_1 \to M$.
We let $M_1'$ denote the manifold obtained by splitting $M'$ along
$S'$, and let $q_1\co M_1'  \to M_1$ denote the quotient map.  We set
$S_1' = q_1^{-1}(S')$.

First suppose that $F$ is a strict essential surface.  We will assume
that $S$ has been chosen among all surfaces isotopic to $F$ to be
transverse to $T$ and to meet $T$ in the minimal number of simple
closed curves.  We will show that $S$ satisfies conditions (1) and
(2).

If $S\cap T=\emptyset$ then $S$ is a strict essential surface in $M'$
by \ref{meetsall}, so conditions (1) and (2) hold in this case.  Thus
we may assume that $S\cap T\not=\emptyset$.  In particular, since
$S$ is connected, no component of $S'$ is closed.  No
component of $S'$ can be boundary-parallel since otherwise
$S$ would be isotopic to a surface that meets $T$ in fewer simple
closed curves.  Since $M$ is irreducible it follows in particular that
no component of $S'$ is a disk.  This, together with the
$\pi_1$--injectivity of $S$, implies that $S'$ is
$\pi_1$--injective in $M'$.  Hence condition (1) holds for $S$.

To prove that condition (2) holds, assume that every component of $S'$
is a semifiber in the component of $M'$ containing it.  We will show
that $M_1$ is a $[0,1]$--bundle, and hence that $S$ is a semifiber,
contradicting our supposition that $S$ is strict.  By \ref{meetsall}
our assumption implies that $S'\cap T'\not=\emptyset$, and hence every
component of $M'$ contains a component of $S'$.  The manifold obtained
by splitting a component $X$ of $M'$ along a component of $S'$ is a
$[0,1]$--bundle.  According to \ref{strictlyvertical}, all of the other
components of $S'$ in $X$ are isotopic to horizontal surfaces in this
$[0,1]$--bundle.  It then follows that the manifold $M_1'$ has the
structure of a $[0,1]$--bundle for which the associated
$\{0,1\}$--bundle is the surface $S_1'$.  The vertical boundary of the
$[0,1]$--bundle $M_1'$ is $V = q_1^{-1}(\partial M')$.  Observe that
$M_1$ is a quotient of $M_1'$ obtained by identifying pairs of
components of $V$.  To obtain the required $[0,1]$--bundle structure on
$M_1$ it suffices to observe that the gluing maps are isotopic to
fiber-preserving maps with respect to the $[0,1]$--bundle structures on
the components of $V$.  This is because any homeomorphism between two
trivial $[0,1]$--bundles over $S^1$ is isotopic to a fiber-preserving
map.

As a preliminary to proving the converse we observe that, since $T$ is
$\pi_1$--injective and $M$ is irreducible, any properly embedded disk
in $M'$ having its boundary contained in $T'$ must be
boundary-parallel in $M'$.

For the proof of the converse we assume that conditions (1) and (2)
hold for the surface $S$.  No component of $S'$ can be a $2$--sphere.
Hence if $S$ had a $2$--sphere component then some component of $S'$
would be a disk whose boundary is contained in $T'$.  Since any
such disk is boundary-parallel, this would contradict condition (1).
Thus $S$ has no $2$--sphere components.

Suppose that $S$ fails to be $\pi_1$--injective.  Then we may choose
a compressing disk $D$ for $S$ which is transverse to $T$ and meets
$T$ in the minimal number of components.  Since $T$ is
$\pi_1$--injective, the minimality implies that all components of
$D\cap T$ must be arcs.  Since $S'$ is $\pi_1$--injective in
$M'$ according to (1), we have $D\cap T \neq \emptyset$.  An outermost arc $\beta$ of $D\cap T$ is the image under
$q$ of an arc in $\partial M'$ which is fixed-endpoint homotopic in
$M'$ to an arc $\alpha$ in $S'$.  It follows from an
application of \fullref{boundaryincompressible}, with $M$ replaced
by $M'$ and $F$ by $S'$, that $\alpha$ is the frontier of a
disk $E_1$ in $S'$.  The subdisk of $D$ cobounded by $\alpha$
and $\beta$ is the image of a disk $E_2$ in $M'$.  The union of $E_1$
and $E_2$ is a properly embedded disk in $M'$ having its boundary
contained in $T'$, and must therefore be a boundary-parallel
disk in $M'$.  It now follows that $D$ is isotopic to a disk that
has fewer components of intersection with $T$, contradicting our
choice of $D$.  This shows that $S$ is $\pi_1$--injective.

Suppose that $S$ is boundary-parallel.  We will show that some
component of $S'$ is boundary-parallel, contradicting condition (1).
There is a submanifold $P$ of $M$ whose frontier is $S$ such that the
pair $(P,S)$ is homeomorphic to $(S\times I, S\times\{1\})$.  If
$T\cap P = \emptyset$ then it is immediate that $S'$ is
boundary-parallel in $M'$.  Since $S$ is $\pi_1$--injective, if some
component of $T\cap S$ is homotopically trivial in $M$ then it must
bound a disk in $S$.  A minimal disk in $S$ bounded by a component of
$T\cap S$ is the image under $q$ of a properly embedded disk $D\subset
M'$ with $\partial D\subset T'$.  Since $D$ must be boundary-parallel
in $M'$, it is the required boundary-parallel component of $S'$.  If
every component of $T\cap S = \partial (T\cap P)$ is homotopically
nontrivial in $M$ then since $T$ is $\pi_1$--injective, $T\cap P$ is
$\pi_1$--injective.  It therefore follows from \ref{strictlyvertical}
that every component of $T\cap P$ is parallel to a subsurface of $S$.
This implies that some component of $S'$ is boundary-parallel in $M'$.
This completes the proof that $S$ is essential.

To show that $S$ is strict assume, to the contrary, that the manifold
$M_1$ is a $[0,1]$--bundle and that the associated $\{0,1\}$--bundle is
$S_1 = p^{-1}(S)$. Then $T_1 = p^{-1}(T)$ is a disjoint union of
annuli in $M_1$ whose boundary components are contained in the
$\{0,1\}$--bundle $S_1$.  If any of these annuli were parallel into
$S_1$ it would imply that some component of $S'$ is a
boundary-parallel annulus in $M'$, contradicting condition (1).  It
therefore follows from \ref{strictlyvertical} that each component of
$T_1$ is isotopic to a vertical annulus in the $[0,1]$--bundle $M_1$.
Therefore the manifold $M_1'$, which we can think of as being obtained
by splitting $M_1$ along $T_1$, is a $[0,1]$--bundle whose associated
$\{0,1\}$--bundle is $S_1'$.  This contradicts condition (2).  Thus
$S$ is a strict essential surface in $M$.  \EndProof

\section{Cable spaces}

We have said that \fullref{mfldwithtorus} can be used to
identify strict essential knots in cabled knot exteriors. This is
because the exterior of a cabling $K'$ of a knot $K$ can be regarded
as being constructed from $M(`K)$ by gluing a ``cable space'' to its
boundary. We begin this section with a definition and some
observations related to this construction, and then prove
\fullref{cablelemma}, which provides a wealth of essential surfaces
in a cable space.

\Number\label{cabledef}
We define a \textit{cable space\/} to be a Seifert fibered manifold over an
annulus with one singular fiber.  Note that a cable space has exactly
three isotopy classes of essential vertical annuli; one has a
boundary curve on each boundary torus of the cable space and the
other two have both boundary curves on the same boundary torus.
\EndNumber

\Number\label{CableKnotCableSpace}   Let $K$ be a knot in a closed,
orientable, connected $3$--manifold and $K'$ be a $q$--strand
cabling of $K$ for some $q\ge2$. It follows from the definition of a cable
in \ref{cable} that $K$ has a regular neighborhood $V$ such
that $K'$ is contained in a torus $T\subset\inter V$ which is
boundary-parallel in $V$, and $K'$ has geometric intersection number $q$
with the boundary of a meridian disk for the
solid torus $V$. We shall call a neighborhood $V$ with these
properties an \textit{enveloping solid torus\/} for the cabling $K'$ of
$K$. 

If $V$ is an enveloping solid torus for a $q$--strand cabling $K'$ of a
knot $K$ then $V$ admits a Seifert fibration over a disk in which $K'$
is a regular fiber, $K$ is the only singular fiber, and the order of
this singular fiber is $q$. Hence if $W\subset\inter V$ is a regular
neighborhood of $K'$ disjoint from $K$, then $N=\overline{V-W}$ admits
a Seifert fibration over an annulus which has exactly one singular
fiber, and the order of the singular fiber is $q$. In particular, $N$
is a cable space.
  \EndNumber

\Lemma\label{cablelemma} Suppose that $N$ is a cable space (see
\ref{cabledef}) with boundary tori $T_1$ and $T_2$.  Let $\iota^j$
denote the inclusion map from $T_j$ to $N$.  Then there exists a bijection
$\phi=\phi_{N,T_1,T_2}\co {\mathcal S}(T_1) \to{\mathcal S}(T_2)$ having the
following properties:

\begin{enumerate}
\item[(i)] if $\alpha_j$ is a primitive element of 
$H_1(T_j;\Z)\subset H_1(T_j;\Q)$ for $j=1,2$, and if
$\phi(\langle\alpha_1\rangle)= \langle\alpha_2\rangle$, then $\iota^2_\sharp(\alpha_2)$ is
a rational multiple of $\iota^1_\sharp(\alpha_1)$ in $H_1(N;\Q)$;

\item[(ii)] for each
slope $\sigma$ on $T_1$ there exists a connected essential surface in $N$
which has nonempty intersection with both $T_1$ and $T_2$ and has $\sigma$
and $\phi(\sigma)$ as boundary slopes.
\end{enumerate}
\EndLemma

\Proof  We
identify $H_1(\partial N;\Q)$ with $H_1(T_1;\Q)\oplus
H_1(T_2;\Q)$.  The cable space $N$ may be decomposed as the union
of homeomorphic copies of $D^2\times S^1$ and $S^1\times S^1\times I$,
meeting along an annulus.  A Mayer--Vietoris computation shows that
${\iota^j}_\sharp$ is an
isomorphism $H_1(T_j;\Q) \to H_1(N;\Q)$ for $j = 1,2$.  If $\alpha_1$ is a primitive element of
$H_1(T_1;\Z)\subset H_1(T_1;\Q)$ then up to sign there is a
unique primitive element $\alpha_2$ of $H_1(T_2;\Z)\subset
H_1(T_2;\Q) $ which is a rational multiple of
$(\iota^2_\sharp)^{-1}\circ\iota^1_\sharp(\alpha_1)$. We then define the map
$\phi$ by setting $\phi(\langle\alpha_1\rangle) = \langle\alpha_2\rangle$; this is a well-defined bijection since every slope
on $T_j$ may be written in the form $\langle\alpha\rangle$, where
$\alpha$ is determined  up to sign by the slope. Property (i) of
$\phi$ is immediate from the definition.

To prove Property (ii), let $\sigma$ be any slope on $T_1$.  Write $\sigma@{=}@\langle\alpha_1\rangle$ and $\phi(\sigma)@{=}@\langle\alpha_2\rangle$ where $\alpha_1 \in H_1(T_1;\Z)$ and
$\alpha_2 \in H_1(T_2;\Z)$ are primitive elements such that the sum
$m\,\iota^1_\sharp(\alpha_1)+n\,\iota^2_\sharp(\alpha_2)=0$ for some
relatively prime integers $m$ and $n$.  By the long exact homology
sequence of $(N,\partial N)$ there is a class $\sea$ in
$H_2(N,\partial N;\Z)$ whose image under the boundary map
$\partial \co  H_2(N,\partial N;\Z) \to H_1(\partial N;\Z)$ is $m\alpha_1
\oplus n\alpha_2$.  It follows from the proof of \cite[Lemma
6.6]{hempel} that there is an oriented essential surface $S$ which
represents the homology class $\sea$.  Since $\partial(\sea) m\alpha_1 \oplus n\alpha_2$, and since $S$ is essential and $T_1$ and
$T_2$ are tori, every component of $\partial S\cap T_1$ has slope
$\sigma$ and every component of $\partial S\cap T_2$ has slope
$\phi(\sigma)$.

Since $\partial(\sea)\not=0$, some component $S_0$ of $S$ must
represent a class $\sea_0 \in H_2(N,\partial N; \Z)$ such that
$\partial(\sea_0)\not= 0$.  Note that $S_0 \cap T_2$ is nonempty
since otherwise $\partial(\sea_0)$ would be a nonzero multiple of
$\alpha_1$, which is impossible since the image of $\alpha_1$ in
$H_1(N;\Z)$ has infinite order.  Similarly $S_0$ must have
nonempty intersection with $T_1$.  Furthermore, since $S_0$ is a
component of $S$ it is essential and has boundary slopes $\sigma$ and
$\phi(\sigma)$.
\EndProof

\section{The effect of cabling on strict boundary slopes}
\label{Asection}

The goal of this section is to prove \fullref{theoremA} of the
Introduction.

\Number\label{cableconventions}
In the next two lemmas, \ref{ExplicitPhi} and \ref{MultipleCopies}, we
shall consider a knot $K$ in a closed, connected, orientable
$3$--manifold $\Sigma$ and   a cabling $K'$ of $K$ with $q\ge2$
strands. We shall denote by $V$
an enveloping solid torus for the cabling $K'$ and by $W$ a
regular neighborhood of $K'$ that is contained in $V$ and disjoint
from $K$. We shall make the explicit choices $V(`K)=V$ and $V(`K')=W$
for the regular neighborhoods of $K$ and $K'$ (see \ref{meridian}). We
shall set $N=\overline{V-W}$, so that $N$ 
is a cable space by \ref{CableKnotCableSpace}.  We set $T_1=\partial
V=\partial M(`K)$ and 
$T_2=\partial W=\partial M(`K')$. Thus $T_1$ and $T_2$ are the boundary
tori of the cable space $N$. We let
$\phi=\phi_{N,T_1,T_2}\co {\mathcal S}(T_1)
\to{\mathcal S}(T_2)$ denote the bijection given by
\fullref{cablelemma}. 
\EndNumber


\begin{lemma}\label{ExplicitPhi}
Suppose that $K$ is a knot in a closed, connected, orientable
$3$--manifold $\Sigma$, and   that $K'$  is a cabling of $K$ with $q\ge2$
strands. Let $V$, $W$, $N$, $T_1$, $T_2$ and $\phi$ be defined as in
\ref{cableconventions}. 
 Let
$\nu=\nu_{\mu,\lambda}\co {\mathcal S}(T_1) \to\Q\cup\{\infty\}$
and $\nu'=\nu_{\mu',\lambda'}\co {\mathcal S}(T_2) \to\Q\cup\{\infty\}$ be the bijections given by \ref{numerical},
 Then there exist constants $u\in\Q$ and
$\epsilon\in\{\pm1\}$ such that for every $\sigma\in{\mathcal S}(T_1)
$ we have
$$\nu'(\phi(\sigma))=\epsilon q^2\nu(\sigma)+u.$$ 
Here the right-hand side is interpreted to be $\infty$ if $\nu(\sigma)=\infty.$
\end{lemma}

\begin{proof}
Since $V$ is an enveloping solid torus for the $q$--strand cabling $K'$ of $K$,
there is a meridian disk $D$ for the solid torus $V$ which
meets $K$ transversally in one point, and meets $K'$ transversally in
$q$ points. Furthermore, the intersections of $D$ with $K'$ are all
 consistently oriented, in the sense that if $\omega$ is a transverse
orientation to $D$ in $V$, the orientations of $K'$ induced by
$\omega$ at the different points of $D\cap K'$ all coincide. By
choosing $D$ to be in standard position with respect to $W$ we may
arrange that $D\cap W$ consists of $q$ meridian disks in $W$,
each of which contains a unique point of $D\cap K'$. Now $P=D\cap
N$ is a planar surface whose boundary consists of one meridian curve
in $T_1=\partial V$ and $q$ meridian curves in $T_2=\partial W$.

We identify $H_1(\partial N;\Z)$ with $H_1(T_1;\Z)\oplus
H_1(T_2;\Z)$. Since the intersections of $D$ with $K'$ are all
 consistently oriented, we
may orient $P$ in such a way that the class $[P]\in H_2(N,\partial
N;\Z)$ satisfies
\begin{equation}
\label{boundaryP}
\partial [P]=\mu+\zeta q\mu'
\end{equation}
 for some
$\zeta\in\{\pm1\}$. 

For $i=1,2$ and for every $\alpha\in H_1(T_i;\Z)$, let us denote by
$\bar\alpha$ the image of $\alpha$ under the natural homomorphism
$H_1(\partial N;\Z)\to H_1(N;\Q)$. 
It follows from \eqref{boundaryP} that
\begin{equation}\label{barmu}
\bar\mu=-\zeta q\bar\mu'.
\end{equation}

We let $\theta\in\{\pm1\}$ denote the homological intersection number
  of  $\lambda'$  with $\mu'$ in $T_2$. Then by \eqref{boundaryP},
  the  homological intersection number of $\bar\lambda'$ with
  $[P]$ in $N$ equals $\zeta\theta
  q$. 

We saw in the proof of \fullref{cablelemma} that the inclusion
homomorphism from $H_1(T_i;\Q)$ into $H_1(N;\Q)$ is an isomorphism
for $i=1,2$. In particular, $\bar\lambda$ and $\bar\mu$ form a basis for
$H_1(N;\Q)$. Let us write 
$\bar\lambda' =t \bar\mu + h \bar\lambda$ for some  $t,h\in\Q$.
If $\eta\in\{\pm1\}$ denotes the homological intersection number
  of  $\lambda$  with $\mu$ in $T_1$, it follows again from
  \eqref{boundaryP}
  that the homological intersection number
  of  $\bar\lambda'$  with $[P]$ in $N$ is $\eta h$. Hence $\eta
  h=\zeta\theta q$, and therefore
\begin{equation}\label{barlambda}
\bar\lambda' =t \bar\mu + \zeta\theta\eta q \bar\lambda.
\end{equation}

We shall show that the lemma holds if we set
$\epsilon=\eta\theta\in\{\pm1\}$ and $u=\zeta qt\in\Q$.

Consider any element $\sigma$ of ${\mathcal S}(T_1)$. We may write
$\sigma=\langle\alpha\rangle$ for some primitive element
$\alpha=a\mu+b\lambda$ of $H_1(T_1;\Z)$, where $a$ and $b$ are relatively
prime integers. We have $\nu(\sigma)=-a/b$. Likewise, if we write
$\phi(\sigma)=\langle\alpha'\rangle$ for some primitive 
$\alpha'=a'\mu'+b'\lambda'$ in $H_1(T_1;\Z)$, then
 $\nu'(\phi(\sigma))=-a'/b'$.
Using \eqref{barmu} and \eqref{barlambda}, we find that
$$\bar\alpha'=a'\bar\mu'+b'\bar\lambda'=-\frac{\zeta a'\bar\mu}{q}+b'(t
\bar\mu + \zeta\theta\eta q \bar\lambda)=\Big(\frac{b'qt-\zeta a'}{q}\Big)\bar\mu+\zeta\theta\eta b'q \bar\lambda.$$

According to \ref{cablelemma} (i), $\bar\alpha$ is a
rational multiple of $\bar\alpha'$, say $\bar\alpha=r\bar\alpha'$. Hence
$$\bar\alpha=\Big(\frac{b'qt-\zeta a'}{q}\Big)r\bar\mu+\zeta\theta\eta b'qr \bar\lambda.$$
On the other hand, we have 
$$\bar\alpha=a\bar\mu+b\bar\lambda,$$
and since $\bar\lambda$ and $\bar\mu$ form a basis for
$H_1(N;\Q)$, it follows that
$$a=\Big(\frac{b'qt-\zeta a'}{q}\Big)r \quad \text{and} \quad b=\zeta\theta\eta b'qr.$$
$$\nu(\sigma)=-\frac{a}{b}=-\frac{(b'qt-\zeta a')/q}{\zeta\theta\eta b'q}
=\frac{\theta\eta}{q^2}\frac{a'}{b'}-
\frac{\zeta\eta\theta t}{q}=\frac{\theta\eta}{q^2}\nu'(\phi(\sigma))-
\frac{\zeta\eta\theta t}{q}.\leqno{\hbox{Hence}}
$$
This gives the required equality 
$$\nu'(\phi(\sigma))=\theta\eta q^2\nu(\sigma)+\zeta qt=\epsilon
q^2\nu(\sigma)+u.\proved$$
\end{proof}

\begin{lemma}\label{MultipleCopies}
Suppose that $K$ is a knot in a closed, connected, orientable
$3$--manifold $\Sigma$, and   that $K'$  is a cabling of $K$ with $q\ge2$
strands. Let $V$, $W$, $N$, $T_1$, $T_2$ and $\phi$ be defined as in
\ref{cableconventions}. Suppose that $\sigma\in{\mathcal S}(T_1)={\mathcal
  S}(\partial M(`K))$ is a
strict boundary slope for $K$. Then $\phi(\sigma)\in{\mathcal S}(T_2)={\mathcal
  S}(\partial M(`K'))$ is a
strict boundary slope for $K'$.
\end{lemma}

\begin{proof}
  Since $\sigma$ is a strict boundary slope for $K$, there is a
  connected strict essential surface $F\subset M(`K)$ having boundary
  slope $\sigma$. On the other hand, by \fullref{cablelemma}, there
  is a connected essential surface $E\subset N$, having nonempty
  intersection with both $T_1$ and $T_2$ and having $\sigma$ and
  $\phi(\sigma)$ as boundary slopes. Let $m$ and $n$ denote,
  respectively, the numbers of components of $\partial F$ and
  $\partial E$. Let $E^*\subset N$ be an essential surface with $m$
  components, each isotopic to $E$ in $N$, and let $F^*\subset M$ be a
  strict essential surface with $n$ components, each isotopic to $F$
  in $M$. Then $\partial E^*$ and $\partial F^*$ are both
  $1$--manifolds in $T_1$ with slope $\sigma$, and each of them has
  $mn$ components. Hence after varying $E$ within its isotopy class we
  may assume that $\partial E=\partial F$. This means that $F'=E\cup
  F$ is a properly embedded surface in $M(`K')=M(`K)\cup N$ transverse
  to $T_1$.  Since
  $E=F'\cap N$ is an essential surface in $N$, and $F=F'\cap M(`K)$ is a
  strict essential surface in $M(`K')$, it follows from \fullref{mfldwithtorus} that $F'$ is a strict essential surface in
  $M(`K')$. As the boundary slope of $F'$ is clearly equal to
  $\phi(\sigma)\in{\mathcal S}(T_2)={\mathcal S}(\partial M(`K'))$, it follows
  that  $\phi(\sigma)$ is a boundary slope for $K'$.
\end{proof}

\begin{proof}[Proof of \fullref{theoremA}]
    We are given a nonround knot $K$ in a closed, connected, orientable
  $3$--manifold $\Sigma$ and a cabling $K'$ of $K$ with $q\ge2$
  strands such that $K'$ is meridionally small.    We define
  $V$, $W$, $N$, $T_1$, $T_2$ and $\phi$ as in \ref{cableconventions}.
  We also choose framings $(\mu,\lambda)$ and $(\mu',\lambda')$ for
  $K$ and $K'$, and let $\nu=\nu_{\mu,\lambda}\co {\mathcal S}(T_1) \to\Q
  \cup\{\infty\}$ and $\nu'=\nu_{\mu',\lambda'}\co {\mathcal
    S}(T_2)\to\Q\cup\{\infty\}$ be the bijections given by
  \ref{numerical}. We denote by $\epsilon\in\{\pm1\}$ and $u\in\Q$ the constants given by \fullref{ExplicitPhi}.
  
  We must first show that $K$ is meridionally small, ie, that $M(`K)$
  is irreducible and that the meridian slope of $K$ is not a strict
  boundary slope. If $S\subset\inter
  M(`K)$ is a $2$--sphere, then $S$ bounds a ball $B$ in the interior of
  $M(`K')=M(`K)\cup N$. Since $N$ is connected and disjoint from $S$, we
  must have either $N\subset\inter B$ or $N\cap B=\emptyset$. But
  $N\subset\inter B$ would imply $\partial M(`K')\subset B$, a
  contradiction. Hence $N\cap B=\emptyset$ and therefore $B\subset
  M(`K)$. This shows that $M(`K)$ is irreducible.
  
Next note that, according to \fullref{ExplicitPhi}, we have 
$\nu'(\phi(\langle\mu\rangle))=\epsilon
q^2\nu(\langle\mu\rangle)+u$,  where $\nu(\langle\mu\rangle)=\infty$.
Hence $\nu'(\phi(\langle\mu\rangle))=\infty$,
that is, $\phi(\langle\mu\rangle)=\mu'$. 
If $\langle\mu\rangle$ were a strict
  boundary slope for $K$, \fullref{MultipleCopies} would now imply that 
the meridian slope $\langle\mu'\rangle$ is a strict
  boundary slope for $K'$, a contradiction to the hypothesis that $K'$
  is meridionally small.  Hence the meridian slope $\langle\mu\rangle$ is not a strict
  boundary slope for $K$. This
completes the proof that $K$ is meridionally small.

It remains to show that
 $d_{K'} \geq q^2d_K$. By definition of $d_K$, there
 exist strict boundary slopes $\sigma$ and $\tau$ for $K$ such that 
$$\nu(\sigma)-\nu(\tau)=d_K.$$
According to \fullref{MultipleCopies}, the slopes
$\sigma'=\phi(\sigma)$ and $\tau'=\phi(\tau)$ are strict boundary
slopes for $K'$. But from \fullref{ExplicitPhi} we have
$$|\nu'(\sigma')-\nu'(\tau')||(\epsilon q^2\nu(\sigma)+u)-(\epsilon q^2\nu(\tau)+u)|=q^2d_K.$$ 
This shows that $d_{K'}\ge q^2d_K$.
\end{proof}

\section[A lower bound for d_K]{A lower bound for $d_K$}
\label{Bsection}

The goal of this section is to prove \fullref{theoremB} of the Introduction.

\begin{proposition}\label{finite}
Suppose that $K$ is a knot in a closed, connected, orientable
  $3$--manifold $\Sigma$. Then for some integer $n\ge0$ there exist
 knots $K_{0}, K_{1}, \ldots, K_{n}$ in $\Sigma$
such that (i) $K=K_{0}$, (ii) $K_{n}$ is not a cable knot and (iii) for
each $i$ with $0 \leq i \leq n-1$, the knot $K_{i}$ is a cabling of
$K_{i{+}1}$.  
\EndProposition

\Remark\label{IteratedRemark} According to our definition of a cable
knot in \ref{cable}, it may happen that the knot $K_n$ given by
\fullref{finite} is round. If $K_n$ is round then it follows
from the definition that $K$ is a generalized iterated torus knot.
\EndRemark

\begin{proof}[Proof of \fullref{finite}]
  Suppose that $K$ is a knot in a closed, connected, orientable
  $3$--manifold $\Sigma$. Set $M=M(`K)$.  Since $M$ is a compact,
  irreducible, orientable $3$--manifold, it follows from Haken's
  finiteness theorem \cite[Lemma 1.32]{hempel}  that there is
  an integer $\Theta>0$ with the following property: if $F_1,\ldots,F_\Theta$
  are disjoint, closed, connected, orientable surfaces of strictly
  positive genus in $\inter M$ such that the homomorphism $\pi_1(F_i)\to\pi_1(M)$ is
  injective for $i=1,\ldots,\Theta$, then the closure of some component of
  $M-(F_1\cup\cdots\cup F_\Theta)$ is homeomorphic to $F\times[0,1]$ for
  some closed surface $F$.
  
  Now assume that the conclusion of \fullref{finite} does not
  hold. We shall recursively construct knots $K_0,\ldots,K_\Theta$ in
  $\Sigma$ and regular neighborhoods $V_i$ of the $K_i$ in such a
  way that for $i=0,\ldots,\Theta{-}1$ these conditions hold: (i) $K_i$
  is a cabling of $K_{i{+}1}$, (ii) $V_{i{+}1}$ is an enveloping solid
  torus (in the sense of \ref{CableKnotCableSpace}) for the cabling
  $K_i$ of $K_{i{+}1}$ and (iii) $V_i$ is contained in the interior of
  $V_{i{+}1}$ and is disjoint from $K_{i{+}1}$.

  We set $K_0=K$ and set $V_0=V(`K)$. Now suppose
    that for a given $m\in\{0,\ldots,\Theta{-}1\}$ we have defined
    $K_0,\ldots,K_m$ and $V_0,\ldots,V_m$ so that
    conditions (i)--(iii) hold for every $i$ with $0\le i<m$. (This
    is of course vacuously true when $m=0$.) We need to define
    $K_{m+1}$ and $V_{m+1}$ so that (i)--(iii) hold for
    $i=m$.

If $K_m$ were not a cable knot, then since condition (i) holds for
$0\le i<m$, the conclusion of \fullref{finite} would hold with
$n=m$. As we have assumed that this conclusion does not hold, $K_m$ is
a cable knot. In particular, $K_m$ is a cabling of some knot
$K^*$ in $\Sigma$. Let $V^*$ be an enveloping solid torus
for the cabling $K_m$ of $K^*$, and let $W^*$ be a regular neighborhood
of $K_m$ which is contained in the interior of
  $V^*$ and is disjoint from $K_{m+1}$. Since $W^*$ and $V_m$ are both
  regular neighborhoods of $K_m$ in $\Sigma$, there is a homeomorphism
  $h\co \Sigma \to\Sigma$, isotopic to the identity by an isotopy fixing
  $K_m$, such that $h(W^*)=V_m$. The knot $K_{m+1}=h(`K^*)$ and the
  solid torus $V_{m+1}=h(V^*)$ then have the required properties. This
  completes the recursive construction.

If $K_\Theta$ were not a cable knot, then since condition (i) holds for
$0\le i<\Theta$, the conclusion of \fullref{finite} would hold with
$n=\Theta$. As we have assumed that this conclusion does not hold, $K_\Theta$ is
a cable knot. In particular, $K_\Theta$ is not round.

For $i=0,\ldots,\Theta{-}1$ set $N_i=\overline{V_{i{+}1}-V_i}$. According to
\ref{CableKnotCableSpace}, each $N_i$ is a cable space.

For $i=0,\ldots,\Theta$, set $T_i=\partial V_i$.  Note that the $T_i$
are disjoint tori contained in $\inter M(`K)$. Set ${\mathcal
  T}=T_1\cup\cdots\cup T_\Theta$. Note that the closures of the
components of $M(`K)-{\mathcal T}$ are $N_0,\ldots,N_{\Theta{-}1}$ and
$\overline{\Sigma-V_\Theta}$. 

We distinguish two cases. First suppose that for some
$i\in\{0,\ldots,\Theta\}$ the inclusion homomorphism
$\pi_1(T_i)\to\pi_1(\overline{\Sigma-V_0})$ has a nontrivial kernel.
Then by 
\cite[Lemma 6.1]{hempel}, there is a disk $D\subset M(`K)$ such that $D\cap{\mathcal
  T}=\partial D$, and such that $\partial D$ does not bound a disk in
$\mathcal T$. Let $Z$ denote the component of $M(`K)-{\mathcal T}$ that contains
$\inter D$, so that $D$ is an essential properly embedded disk in
$\bar Z$. If $\bar Z=\overline{\Sigma-V_\Theta}$, then
$M(`K_\Theta)\cong\overline{\Sigma-V_\Theta}$ contains an essential
disk and therefore has a solid torus as a connected summand; this
contradicts the fact that $K_\Theta$ is not round. If $\bar Z=N_i$ for
some $i\in\{0,\Theta{-}1\}$ we again obtain a contradiction, because the
cable space $N_i$ is a Seifert fibered space with two boundary
components, and the only Seifert fibered space that contains an
essential disk is the solid torus.

There remains the case that 
$\pi_1(T_i)\to\pi_1(\overline{\Sigma-V_0})$ is injective for each
$i\in\{0,\ldots,\Theta\}$.  The defining property of $\Theta$ then
implies that the closure of some component $Z$ of $M(`K)-{\mathcal T}$ is homeomorphic to
$T^2\times[0,1]$. We cannot have $\bar
Z=\overline{\Sigma-V_\Theta}$, since $\overline{\Sigma-V_\Theta}$ has
connected boundary. Hence we must have $\bar
Z=N_i$ for some $i\in\{0,\Theta{-}1\}$. But since the cable space $N_i$ is a
Seifert fibered space over an annulus with a singular fiber, the
fundamental group of either component of $\partial N_i$ is mapped by
the inclusion homomorphism onto a proper subgroup of $\pi_1(N_i)$;
hence $N_i$ cannot be homeomorphic to $T^2\times[0,1]$, and we have a
contradiction in this case as well.
\EndProof

\begin{proof}[Proof of \fullref{theoremB}]
  Given a closed, connected, orientable $3$--manifold $\Sigma$
  such that $\pi_1(\Sigma)$ is cyclic, and a meridionally small knot $K\subset\Sigma$, we must show that either (i) $d_{K} \geq 2$ or
  (ii) $K$ is a generalized iterated torus knot.

Let $n\ge0$ be the integer and $K_{0}, K_{1}, \ldots, K_{n}$ the
knots given by \fullref{finite}. For each $i$, $0 \leq i
\leq n-1$, the knot $K_{i}$ is a $q_{i{+}1}$--strand cabling of $K_{i{+}1}$
for some $q_i\ge2$.

If $K_m$ is round for some $m\le n$, then it follows from \fullref{IteratedRemark} that $K$ is a
generalized iterated torus knot. Thus conclusion (ii) holds in this
case.

Now suppose that none of the knots $K_0,\ldots,K_n$ is round. By $n$
successive applications of \fullref{theoremA} we see that each of the 
$K_i$ is  meridionally small, and that for $i=0,\ldots,n-1$. 
Hence
$$d_K=d_{K_0}\ge q_1^2\cdots q_n^2d_{K_n}.$$

On the other hand, since $\pi_1(\Sigma)$ is cyclic, and since $K_n$ is
meridionally small and is not a round knot or a cable knot, it follows
from Theorem 1.1 of \cite{slopediff} that $d_{K_n}\ge2$.  Hence
$$d_K\ge 2q_1^2\cdots q_n^2\ge2,$$ and so
(i) holds.
\EndProof

\bibliographystyle{gtart}
\bibliography{link}

\end{document}